\documentclass[11pt]{amsart}

\usepackage{amsmath, amssymb, hyperref, amsthm, graphicx, amsaddr}

\title{The entropy of the Angenent torus is approximately $1.85122$}
\author{Yakov Berchenko-Kogan}
\address[Y. Berchenko-Kogan]{Washington University in St.~Louis}
\subjclass[2010]{53C44, 70H25, 65-05, 65Pxx}
\keywords{mean curvature flow, Angenent torus, discrete Euler--Lagrange equations, variational integrator}

\newcommand{\abs}[1]{\left\lvert{#1}\right\rvert}
\newcommand{\norm}[1]{\left\lVert{#1}\right\rVert}
\newcommand{\dd}[2][]{\frac{d{#1}}{d{#2}}}
\newcommand{\pp}[2][]{\frac{\partial{#1}}{\partial{#2}}}

\newtheorem{theorem}{Theorem}[section]

\theoremstyle{definition}
\newtheorem{definition}[theorem]{Definition}

\theoremstyle{remark}
\newtheorem{example}[theorem]{Example}
\newtheorem{question}[theorem]{Question}

\begin{document}
\maketitle

\begin{abstract}
  To study the singularities that appear in mean curvature flow, one must understand \emph{self-shrinkers}, surfaces that shrink by dilations under mean curvature flow. The simplest examples of self-shrinkers are spheres and cylinders. In 1989, Angenent constructed the first nontrivial example of a self-shrinker, a torus. A key quantity in the study of the formation of singularities is the \emph{entropy}, defined by Colding and Minicozzi based on work of Huisken. The values of the entropy of spheres and cylinders have explicit formulas, but there is no known formula for the entropy of the Angenent torus. In this work, we numerically estimate the entropy of the Angenent torus using the discrete Euler--Lagrange equations.
\end{abstract}

\section{Introduction}\label{sec:intro}
Mean curvature flow is a well-studied geometric flow; see the survey paper \cite{cmp15}. Under mean curvature flow, a hypersurface $\Sigma\subset\mathbb R^n$ evolves in time in such a way as to decrease its area as quickly as possible. Some surfaces, such as spheres or cylinders, evolve under mean curvature flow by dilation; such surfaces are known as \emph{self-shrinkers}. Self-shrinkers are key to understanding the singularities that develop under mean curvature flow: Just before a singularity forms, the surface resembles a self-shrinker near the singular point. Colding and Minicozzi proved in \cite{cm12} that for a generic singularity, the corresponding self-shrinker will be a sphere or a cylinder. However, other self-shrinkers are possible. Angenent proved the existence of the first nontrivial self-shrinker in \cite{a92}; it has since been named the \emph{Angenent torus}. While Colding and Minicozzi's result prohibits the Angenent torus from appearing in the limit for a generic singularity, for a generic finite-dimensional family of singularities, it is expected that nontrivial self-shrinkers such as the Angenent torus may appear.

To prove their result, Colding and Minicozzi defined a quantity called the \emph{entropy} of a surface, based on Huisken's $F$-functional \cite{h90}. For a self-shrinker, entropy remains constant during mean curvature flow; for any other surface, the entropy is monotonically decreasing. As such, the initial entropy of a surface limits the kinds of self-shrinkers that may appear at singularities as the surface evolves under mean curvature flow: such a limiting self-shrinker must have smaller entropy than that of the initial surface.

Self-shrinkers are critical points for the entropy. For a self-shrinker that will collapse to the origin after one unit of time, Colding and Minicozzi's entropy coincides with Huisken's $F$-functional, which for surfaces $\Sigma\subset\mathbb R^3$ is the weighted surface area
\begin{equation*}
  \frac1{4\pi}\int_{\Sigma}e^{-\abs x^2/4}\,d\text{Area}.
\end{equation*}
The Angenent torus is thus seen as a critical surface of this weighted area functional. The Angenent torus is rotationally symmetric about the $z$-axis, so we can understand it by understanding its cross-section, a closed curve $\gamma$ in the $(r,z)$-plane. After reducing by the rotational symmetry, the $F$-functional becomes the weighted length
\begin{equation*}
  \frac12\int_\gamma re^{-(r^2+z^2)/4}\,d\text{Arclength}.
\end{equation*}
The cross-section of a rotationally symmetric self-shrinker is a critical curve for this weighted length funcional, or, equivalently, a geodesic with respect to the Riemannian metric
\begin{equation}\label{eq:g}
  g = \frac14\left(r^2e^{-(r^2+z^2)/2}\right)(dr^2+dz^2).
\end{equation}
Angenent proved the existence of his self-shrinking torus by showing that the plane with this metric had a closed geodesic. Computing the length of this geodesic gives the entropy of the Angenent torus, which, as described above, places a lower bound on the entropy of an initial surface that can develop an Angenent torus singularity.

To numerically compute the location of this geodesic in the plane in a way that will let us accurately measure its length, we use Lagrangian mechanics. The intuition for using mechanics to attack this problem is that if a particle is given some initial velocity, it will travel along a geodesic with constant speed. The length of the geodesic is simply the product of the particle's speed and its travel time. With this approach, instead of seeking unparametrized closed curves that are critical points of length, we seek parametrized closed curves $q(t)$, $0\le t\le T$, with velocity $\dot q(t)$, that are critical points for the \emph{action} $\mathfrak S(q)=\int_0^T\frac12\abs{\dot q}^2\,dt$. Such critical curves will be geodesics parametrized by a constant multiple of arclength. %Intuitively, we now seek critical points for the length of the curve squared.

In the 1960s and 1970s, numerical analysts developed \emph{variational integrators}, a new class of methods for numerically computing trajectories in Lagrangian mechanics, with excellent theoretical properties, such as preservation of conserved quantities. See Marsden and West \cite{mw01} for a complete treatment. The key idea of this method is to discretize the action. Instead of seeking closed curves $q(t)$ that are critical points for the action $\mathfrak S(q)$, we instead seek discrete curves, that is, sequences of points $q_0,q_1,\dotsc,q_{N-1},q_N=q_0$ that are critical points for a \emph{discrete action} $\mathfrak S_d(q_1,\dotsc,q_N)$. For a particular choice of $\mathfrak S_d$, called the \emph{exact discrete action} and denoted $\mathfrak S_d^E$, a critical sequence $q_k$ for $\mathfrak S_d^E$ coincides with a critical curve $q$ for $\mathfrak S$ in the sense that the $q_k$ are equally spaced points on the curve $q$, and $\mathfrak S_d^E(q_1,\dotsc,q_N)=\mathfrak S(q)$. In practice, however, the exact discrete action cannot be computed numerically, so we instead work with a discrete action $\mathfrak S_d$ that approximates the exact discrete action $\mathfrak S_d^E$. We compute the critical discrete trajectory for $\mathfrak S_d$, giving us an approximation for the exact trajectory and the exact length.

\begin{figure}
  \centering
  \includegraphics[scale=.8]{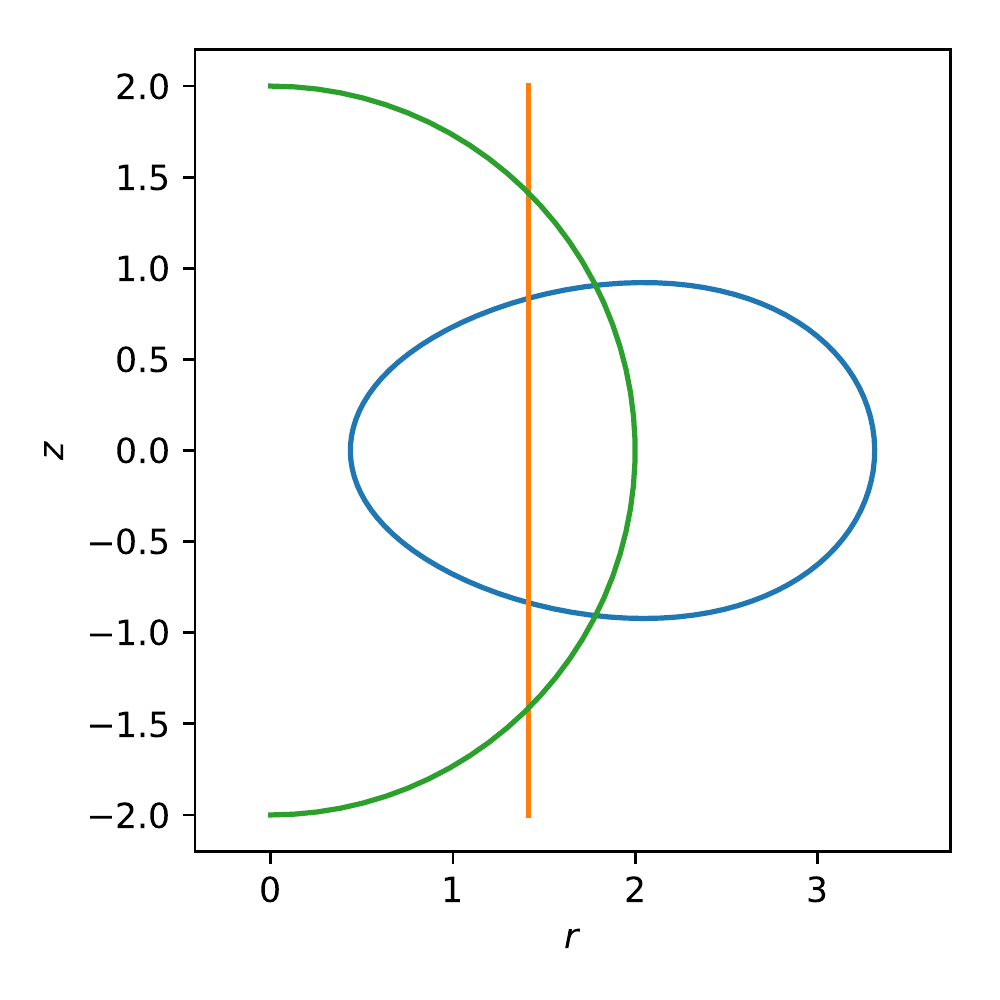}
  \caption{Cross-sections of three self-shrinkers that become extinct in one unit of time: the Angenent torus, the cylinder, and the sphere.}
  \label{fig:angenentcylindersphere}
\end{figure}

In Figure \ref{fig:angenentcylindersphere}, we present the cross-section of the Angenent torus that we computed numerically with $2048$ points, with the cylinder and sphere for comparison. The entropy of the Angenent torus as estimated from the length of this curve is $1.8512186$, and numerical convergence evidence suggests that this value overestimates the true entropy by about $0.0000019$.

Our results are consistent with work of Drugan and Nguyen \cite{dn18}, who proved that the entropy of the Angenent torus is less than $2$, and with the work of Chopp \cite{c94}, who numerically computed a graph of the Angenent torus cross-section using level set methods but did not provide its length. There are many directions in which to continue this work; see Section \ref{sec:futurework}.

\section{Motivation: mean curvature flow}

Mean curvature flow describes the evolution of a surface in an ambient manifold via the negative gradient flow of area. It describes, for example, the evolution of soap films. For the purposes of this exposition, we will restrict our attention to hypersurfaces $\Sigma\subset\mathbb R^{n+1}$. See the survey paper \cite{cmp15} for a more detailed treatment.%If the surface locally resembles a sphere, it will evolve so as to shrink the sphere. If the surface locally resembles a saddle, it will evolve in the direction where the saddle curves more steeply.

\begin{definition}
  The \emph{mean curvature} $H$ at a point $x$ on a surface $\Sigma\subset\mathbb R^{n+1}$ is the trace of the second fundamental form at $x$. Equivalently, the mean curvature is the sum of the principal curvatures of $\Sigma$ at $x$.
\end{definition}
Note that in this convention we do not divide the sum of the principal curvatures by $n$ when computing the mean curvature, so it is more a curvature sum than a curvature mean. The sign convention for $H$ is such that the mean curvature of the sphere is positive.

\begin{definition}
  A family of surfaces $\Sigma_t\subset\mathbb R^{n+1}$ is said to evolve by \emph{mean curvature flow} if the points $x(t)\in\Sigma_t$ move with speed $H$ in the inward normal direction.
\end{definition}

For closed surfaces, this flow will develop a singularity in finite time. For example, a round sphere of radius $r$ will shrink to a point in $\frac{r^2}{2n}$ units of time. The sphere is an example of a surface which evolves under mean curvature flow by dilations. Informally, one could call any such a surface a self-shrinker, but in the literature and in this article we fix the location of the singularity in space and time and use the term \emph{self-shrinker} to refer to a surface that evolves under mean curvature by dilations about the origin and shrinks to a point after one unit of time.

\begin{definition}
  A family of surfaces $\Sigma_t\subset\mathbb R^{n+1}$ for $t<0$ evolving under mean curvature flow is a \emph{self-shrinker} if
  \begin{equation*}
    \Sigma_t=\sqrt{-t}\Sigma_{-1}.
  \end{equation*}
  We will also use the terminology \emph{self-shrinker} for the surface $\Sigma_{-1}$ at time $-1$.
\end{definition}

%Note that, in this definition, we fix the singularity of the self-shrinker to occur one unit of time in the future at the origin.

\begin{example}
  The $n$-dimensional sphere $S^n_{\sqrt{2n}}$ of radius $\sqrt{2n}$ is a self-shrinker, as is the plane $\mathbb R^n$. We also have self-shrinking cylinders $S^k_{\sqrt{2k}}\times\mathbb R^{n-k}$ for every $k$ in between. The earliest nontrivial example of a self-shrinker is a rotationally symmetric $S^1\times S^{n-1}$ proved to exist by Angenent in \cite{a92}. Since then, many other self-shrinkers have been found \cite{cmp15}.

  However, if we restrict our attention to self-shrinkers with $S^{n-1}$ rotational symmetry, Kleene and M\o ller \cite{km14} show that the sphere, the plane, the cylinder $\mathbb R\times S^{n-1}$, and the doughnut $S^1\times S^{n-1}$ are the only complete embedded self-shrinkers, with the caveat that they do not rule out the existence of additional doughnuts other than Angenent's example.
\end{example}

Singularities of mean curvature flow are modeled by self-shrinkers, in the sense that if we take a general surface that develops a singularity as it evolves by mean curvature flow and zoom in at the singular point at just before the singular time, the surface we will see will be close to a self-shrinker. See \cite{cmp15} for a more thorough discussion of this topic.

Huisken \cite{h90} analysed the singularities of mean curvature flow using a monotonicity result for a Gaussian-weighted area called the $F$-functional.

\begin{definition}
  For a surface $\Sigma\subset\mathbb R^{n+1}$, the \emph{$F$-functional} $F(\Sigma)$ is the weighted area
  \begin{equation*}
    F(\Sigma)=(4\pi)^{-n/2}\int_\Sigma e^{-\abs x^2/4}\,d\text{Area}.
  \end{equation*}
  We can, of course, also center the $F$-functional at a different point $x_0$ and time $t_0$ to define
  \begin{equation*}
    F_{x_0,t_0}(\Sigma)=(4\pi t_0)^{-n/2}\int_\Sigma e^{-\abs{x-x_0}^2/(4t_0)}\,d\text{Area}.
  \end{equation*}
  With this notation, $F=F_{0,1}$.
\end{definition}

One motivation for this definition is that the critical points of $F$ are self-shrinkers. That is, if we perturb a self-shrinker, the value of $F$ will remain constant to first order. If we translate and dilate a self-shrinker so that the flow becomes extinct at $x_0$ in time $t_0$, then it will be a critical point for $F_{x_0,t_0}$. The scaling factor $(4\pi)^{-n/2}$ is chosen so that the $F$-functional of a plane through the origin is one.

Colding and Minicozzi \cite{cm12} took this definition one step further and defined the entropy of a surface as the supremum of the $F$-functionals.
\begin{definition}
  The \emph{entropy} of a surface $\Sigma\subset\mathbb R^{n+1}$, denoted $\lambda(\Sigma)$, is the supremum
  \begin{equation*}
    \lambda(\Sigma)=\sup_{x_0,t_0}F_{x_0,t_0}(\Sigma).
  \end{equation*}
\end{definition}

The motivation for this definition is that the entropy is nonincreasing as a general surface evolves under mean curvature flow, and this remains true in the limit as we approach a singularity, in the sense that the self-shrinker that appears when we zoom in at the singular point just before the singular time must have an entropy that is less than or equal to the entropy of the initial surface. Again, see \cite{cmp15} for a more detailed exposition.

For a self-shrinker, $\lambda(\Sigma)=F(\Sigma)$. The Angenent torus, being a self-shrinker, is a critical point for $F$, and our task is to compute the corresponding critical value.

\section{Background: the discrete Euler-Lagrange equations}\label{sec:del}
In Lagrangian mechanics, we consider the motion of a particle in a configuration manifold $Q$. For our purposes, $Q$ will be the half-plane $\{(r,z)\mid r>0\}$. The physics of the system are described by a function $L\colon TQ\to\mathbb R$ called the \emph{Lagrangian}, representing the kinetic energy minus the potential energy. In our case, our Lagrangian is the kinetic energy of a particle of unit mass, with no potential energy term. This Lagrangian is
\begin{equation*}
  L(q,\dot q)=\tfrac12\norm{\dot q}^2_g,
\end{equation*}
where $g$ is the Riemannian metric in \eqref{eq:g}. We then integrate the Lagrangian over a parametrized curve.

\begin{definition}
  For a curve $q(t)$ parametrized by $t\in[0,T]$, the \emph{action} $\mathfrak S(q)$ is defined to be
  \begin{equation*}
    \mathfrak S(q)=\int_0^TL(q,\dot q)\,dt.
  \end{equation*}
\end{definition}

The motivation for this definition is that the curves that describe the motion of a particle in this system are the critical points of $\mathfrak S$ with respect to variations that fix the endpoints. That is, if we perturb the critical curve while leaving the endpoints fixed, the value of $\mathfrak S$ should not change to first order. Such critical curves satisfy a system of differential equations called the Euler--Lagrange equations.
\begin{definition}
  The \emph{Euler--Lagrange equations} are
  \begin{equation}\label{eq:el}
    \dd t\pp[L]{\dot q}=\pp[L]q.
  \end{equation}
\end{definition}

For our choice of Lagrangian, solutions to the Euler--Lagrange equations are geodesics with respect to the Riemannian metric $g$, parametrized with constant speed with respect to $g$-arclength.

We now turn to the question of computing these solutions numerically. While one can certainly apply a standard ODE solver to the Euler--Lagrange equations, there are many advantages to the variational approach, such as preservation of conserved quantities. See the survey paper \cite{mw01} for a more thorough discussion.

The main idea of this approach is to discretize the action. Namely, we split the time interval $[0,T]$ into $N$ small equal time intervals $[t_{k-1},t_k]$ of length $\tau$, and write
\begin{equation*}
  \mathfrak S(q)=\sum_{k=1}^N\int_{t_{k-1}}^{t_k}L(q,\dot q)\,dt.
\end{equation*}
We then approximate $\int_{t_{k-1}}^{t_k}L(q,\dot q)\,dt$ based on the location of the endpoints $q(t_{k-1})$ and $q(t_k)$.

\begin{definition}
  A \emph{discrete Lagrangian} is a function $L_d\colon Q\times Q\to\mathbb R$, chosen so that
  \begin{equation*}
    \int_{t_{k-1}}^{t_k}L(q,\dot q)\,dt\approx L_d\left(q(t_{k-1}), q(t_k)\right).
  \end{equation*}
\end{definition}

In our case, we will make a fairly standard choice of $L_d(q_0,q_1)$. Given a curve $q$ going from $q_0$ to $q_1$ in time $\tau$, we can approximate the position $q$ by the average position $\frac12(q_0+q_1)$, and we can approximate the velocity $\dot q$ by the average velocity $(q_1-q_0)/\tau$, giving us the discrete Lagrangian
\begin{equation}\label{eq:discreteLagrangian}
  L_d(q_0,q_1)=\tau L\left(\frac{q_0+q_1}2,\frac{q_1-q_0}\tau\right).
\end{equation}
We can then write the corresponding action.
\begin{definition}
  For a sequence of points $q_0\dotsc,q_N$, the \emph{discrete action} $\mathfrak S_d$ is the sum
  \begin{equation*}
    \mathfrak S_d(q_0,\dotsc,q_N)=\sum_{k=1}^NL_d(q_{k-1},q_k).
  \end{equation*}
\end{definition}
By design, if $q$ is a curve and we pick $q_k=q(t_k)$ to be equally spaced points on the curve, then $\mathfrak S_d(q_0,\dotsc,q_N)\approx\mathfrak S(q)$. Thinking of the sequence of points $q_0,\dotsc,q_N$ as a discrete curve, we proceed as in the continuous curve setting and seek a discrete curve that is a critical point for $\mathfrak S_d$ with respect to variations that fix the endpoints $q_0$ and $q_N$. A computation shows that such critical discrete curves satisfy a nonlinear system of equations called the discrete Euler--Lagrange equations.
\begin{definition}
  The \emph{discrete Euler--Lagrange equations} are
  \begin{equation}\label{eq:del}
    D_1L_d(q_{k-1},q_k)+D_0L_d(q_k,q_{k+1})=0,\qquad 1\le k\le N-1,
  \end{equation}
  where $D_0L_d$ denotes the derivative of $L_d(q_0,q_1)$ with respect to the $q_0$ variable and $D_1L_d$ denotes the derivative of $L_d$ with respect to the $q_1$ variable. Equivalently, $D_0L_d$ is the gradient of the map $q_0\mapsto L_d(q_0,q_1)$, and $D_1L_d$ is the gradient of the map $q_1\mapsto L_d(q_0,q_1)$.
\end{definition}

Because $\mathfrak S_d$ is an approximation for $\mathfrak S$, solutions to the discrete Euler--Lagrange equations will be approximations for solutions to the Euler--Lagrange equations, in the above sense. There is, however, a particular choice of discrete Lagrangian, called the exact discrete Lagrangian, such that solutions to the discrete Euler--Lagrange equations lie exactly on solutions to the Euler--Lagrange equations.

\begin{definition}\label{def:exactDiscreteLagrangian}
  The \emph{exact discrete Lagrangian} $L_d^E(q_0,q_1)$ for two nearby points $q_0$ and $q_1$ is defined by
  \begin{equation*}
    L_d^E(q_0,q_1)=\int_0^\tau L(q,\dot q)\,dt,
  \end{equation*}
  where $q$ is the solution to the Euler--Lagrange equations going from $q_0$ to $q_1$ in time $\tau$.
  We will call the corresponding discrete action the \emph{exact discrete action} and denote it $\mathfrak S_d^E$.
\end{definition}

For our situation, this curve $q$ is the constant-speed geodesic joining $q_0$ and $q_1$. Its speed is therefore $d(q_0,q_1)/\tau$, where $d$ denotes the distance with respect to the metric $g$. Thus, $L(q,\dot q)=\frac12\left(\frac{d(q_0,q_1)}{\tau}\right)^2$, constant in $t$, and so $L_d^E(q_0,q_1)=\tau L(q,\dot q)=\frac{d(q_0,q_1)^2}{2\tau}$.

Whereas the critical discrete curves for the discrete action $\mathfrak S_d$ lie approximately on the critical curves for the action $\mathfrak S$, the critical discrete curves $\{q_k\}$ for the exact discrete action $\mathfrak S_d^E$ lie \emph{exactly} on the critical curves $q(t)$ for the action $\mathfrak S$, in the sense that $q_k=q(t_k)$. While the exact discrete Lagrangian $L_d^E$ generally cannot be computed numerically, one can estimate the accuracy of the trajectories obtained from a chosen discrete Lagrangian $L_d$ by estimating how close $L_d$ is to $L_d^E$. See \cite[Section 2.3]{mw01}.

\section{Computational methods}
Let $\Sigma\subset\mathbb R^3$ be a self-shrinker that is rotationally symmetric about the $z$-axis. As discussed in Section \ref{sec:intro} and in \cite{a92, dn18}, the cross-section of $\Sigma$ is a geodesic in the $(r,z)$ half-plane with the Riemmanian metric $g$ given by \eqref{eq:g}. Figure \ref{fig:angenentcylindersphere} illustrates three such geodesics, corresponding to the Angenent torus, the cylinder, and the sphere. As discussed in Section \ref{sec:del}, geodesic curves are the solutions to the Euler--Lagrange equations \eqref{eq:el} using the Lagrangian $L(q,\dot q)=\frac12\norm{\dot q}^2_g$. In our numerical analysis, we construct discrete geodesics, that is, sequences of points $q_0,\dotsc,q_N$ in the $(r,z)$ half-plane that approximate the geodesic curves. We do so by solving the discrete Euler--Lagrange equations \eqref{eq:del} with the discrete Lagrangian in \eqref{eq:discreteLagrangian}. Since the path along which a free particle travels does not change if its velocity is scaled, the choice of time step $\tau$ in \eqref{eq:discreteLagrangian} will not affect the discrete trajectories. Thus, for simplicity, we set $\tau=1$.

As with differential equations, to solve the discrete Euler--Lagrange equations, we must specify boundary conditions or initial conditions. We will make use of three such setups.
\begin{itemize}
\item 
  We can fix the endpoints $q_0$ and $q_N$ and solve for the discrete geodesic joining them. That is, we can solve the system of $N-1$ vector equations
  \begin{equation}
    D_1L_d(q_{k-1},q_k)+D_0L_d(q_k,q_{k+1})=0,\qquad 1\le k\le N-1,\label{eq:open}
  \end{equation}
  for the $N-1$ vectors $q_1,\dotsc,q_{N-1}$.
\item 
  We can solve for a closed geodesic by setting $q_0=q_N$, $q_1=q_{N+1}$ and solving the system of $N$ vector equations
  \begin{equation}
    D_1L_d(q_{k-1},q_k)+D_0L_d(q_k,q_{k+1})=0,\qquad 1\le k\le N,\label{eq:closed}
  \end{equation}
  for the $N$ vectors $q_1,\dotsc,q_N$.

  Note that this system of equations must be degenerate because a closed geodesic can always be reparametrized by moving the starting point; that is, solutions are not isolated and instead come in one-parameter families. However, this degeneracy does not appear to cause any difficulties.
  
\item 
  We can specify two initial points, $q_0$ and $q_1$. Then, given $q_{k-1}$ and $q_k$, we can recursively solve
  \begin{equation}
    D_1L_d(q_{k-1},q_k)+D_0L_d(q_k,q_{k+1})=0,\qquad 1\le k\label{eq:shoot}
  \end{equation}
  for $q_{k+1}$. Doing so corresponds to setting an initial position of $(q_0+q_1)/2$ and an initial velocity of $q_1-q_0$, and then computing that geodesic forward in time. Unlike the previous two procedures, this procedure is fast, as it does not require solving a large system of equations.
\end{itemize}

\begin{figure}
  \centering
  \includegraphics[scale=.8]{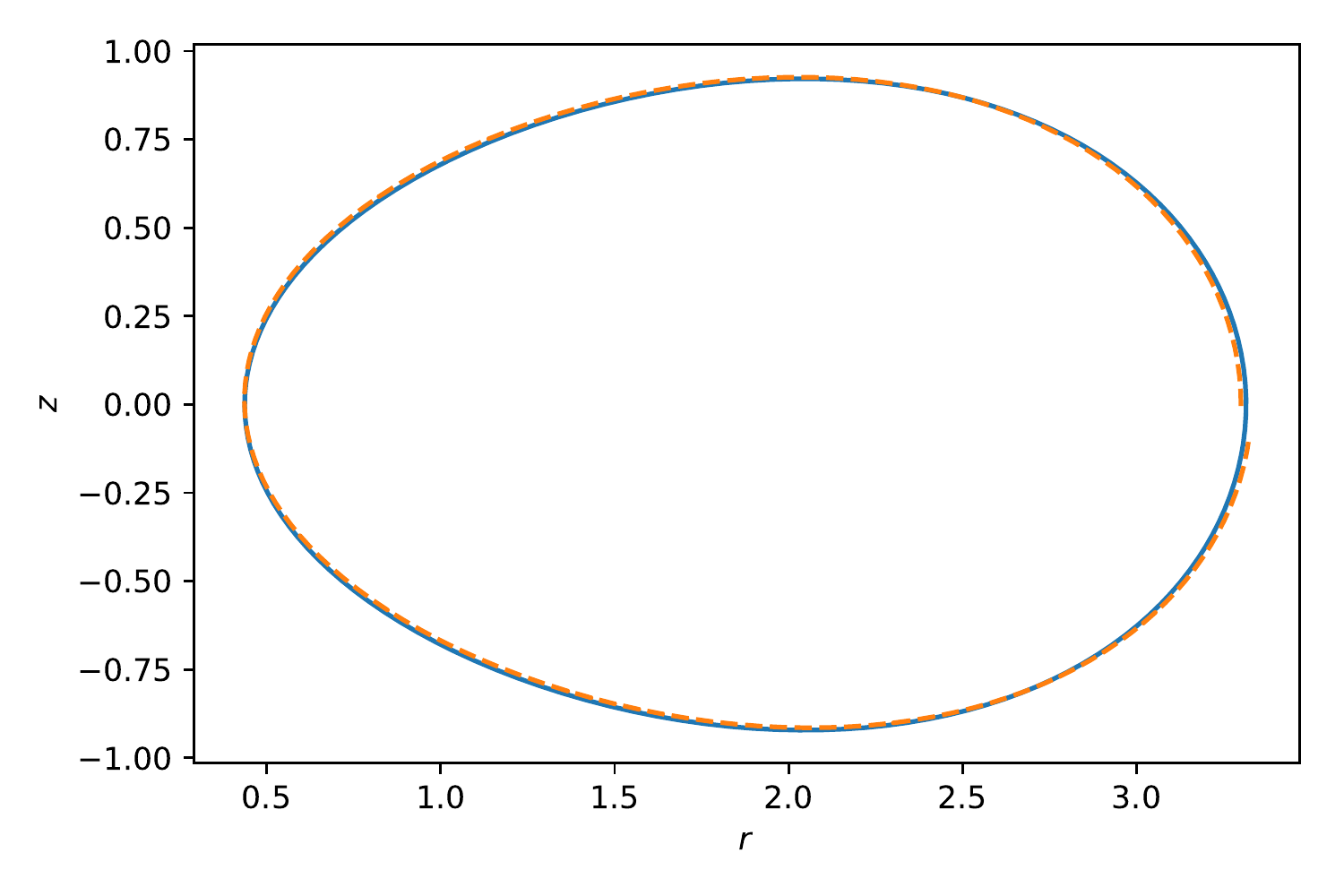}
  \caption{The cross-section of the Angenent torus (solid line), computed using $2048$ points, along with the open geodesic (dashed line) used as the starting estimate for the nonlinear solver. Note that the starting estimate trajectory is off from the correct $r$-intercept, and note that it stops slightly short of making it all the way around.}
  \label{fig:angenent2048withInitialGuess}
\end{figure}

To find the cross-section of the Angenent torus, we solved the system \eqref{eq:closed} for a closed geodesic using SciPy's \texttt{fsolve} routine for solving a nonlinear system of equations. To supply a starting estimate for the closed geodesic required by \texttt{fsolve}, we iteratively solved \eqref{eq:shoot} for an open geodesic with initial data $q_0$ and $q_1$ chosen to be close enough to points on a closed trajectory. Specifically, we chose $q_0=(3.3, -8.5/N)$ and $q_1=(3.3, 8.5/N)$. To find these values, we first experimentally sought a value of $r$ so that a discrete geodesic with initial position at $(r,0)$ and vertical initial velocity appeared to pass by close to its starting point after one time around. Then, we determined the appropriate Euclidean initial speed so that the discrete geodesic completed approximately one cycle after $N$ iterations of \eqref{eq:shoot}. See Figure \ref{fig:angenent2048withInitialGuess} for a plot of the open geodesic used as our starting estimate alongside the closed geodesic found using that starting estimate.

To compute the entropy of the Angenent torus, that is, the length of this geodesic, we estimated the distance $d(q_{k-1},q_k)$ between two consecutive points $q_{k-1}$ and $q_k$. Let $g_{\text{mid}}$ be the value of the metric $g$ at the Euclidean midpoint $(q_{k-1}+q_k)/2$ of these two points. Then we can estimate
\begin{equation*}
  d(q_{k-1},q_k)\approx\norm{q_k-q_{k-1}}_{g_{\text{mid}}}.
\end{equation*}
Of course, $\norm{q_k-q_{k-1}}_{g_{\text{mid}}}$ is precisely the same as $\sqrt{2\tau L_d(q_{k-1},q_k)}$. Indeed, from the discussion following Definition \ref{def:exactDiscreteLagrangian}, we have that
\begin{equation*}
  \frac{d(q_{k-1},q_k)^2}{2\tau}=L_d^E(q_{k-1},q_k)\approx L_d(q_{k-1},q_k).
\end{equation*}
Recall that we can set $\tau=1$ for simplicity. Thus, once we have our discrete closed geodesic $q_0,q_1,\dotsc,q_N=q_0$, we can estimate the entropy of the Angenent torus by computing
\begin{equation*}
  \lambda(\text{Angenent torus})\approx\sum_{k=1}^N\sqrt{2L_d(q_{k-1},q_k)}.
\end{equation*}

\section{Results}
\begin{figure}
  \centering
  \includegraphics[scale=.8]{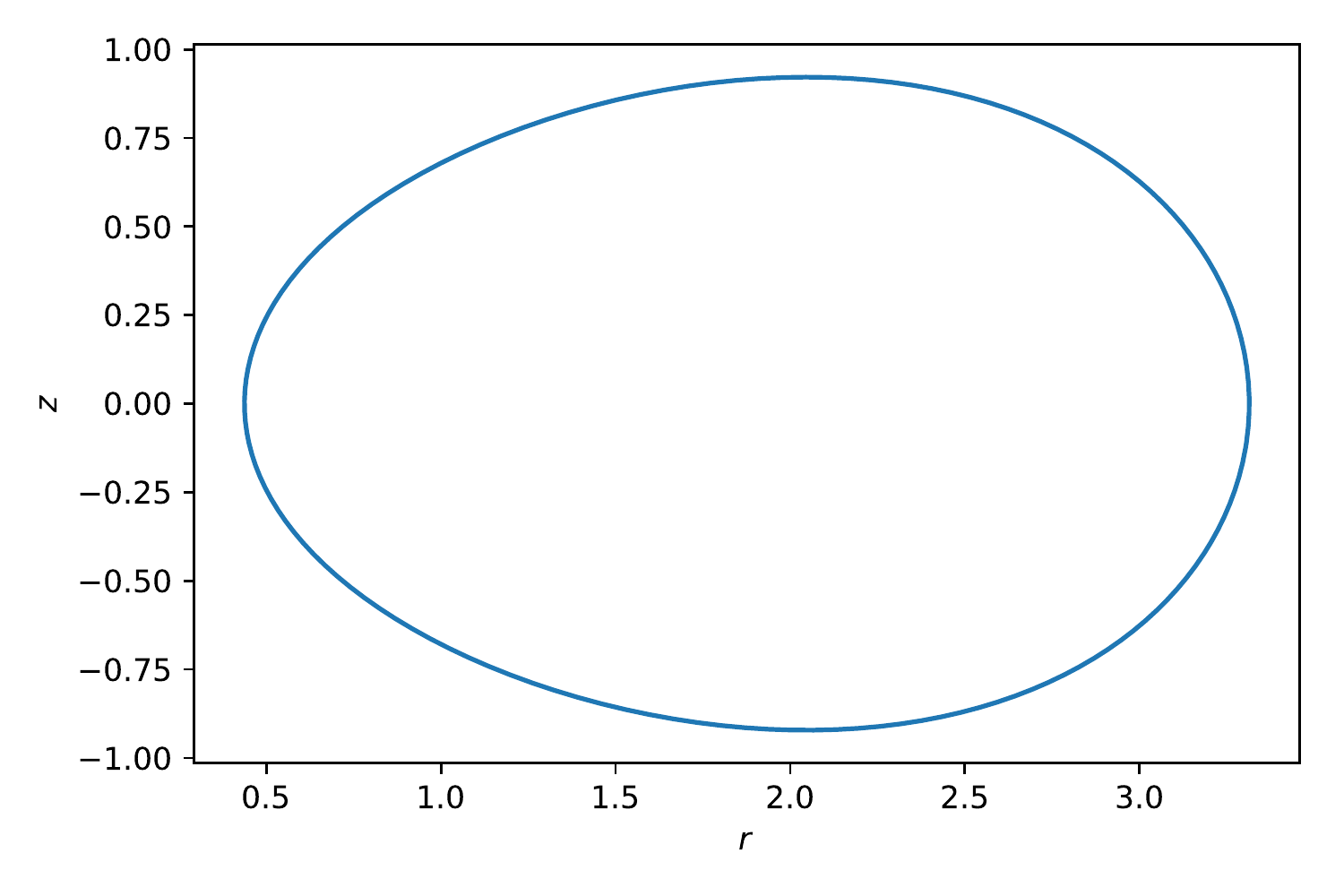}
  \caption{The cross-section of the Angenent torus. The plotted curve was computed using $2048$ points, but, at the resolution of this page, the curves computed using $128$, $256$, $512$, $1024$, and $2048$ points would all appear identical.}
  \label{fig:angenent2048}
\end{figure}

Using $N=128$ points was sufficient for this method to successfully find a closed discrete geodesic. We also performed this computation using $256$, $512$, $1024$, and $2048$ points. See Figure \ref{fig:angenent2048} for the result. The torus crosses $z=0$ at $r=0.4371$ and at $r=3.3147$. It reaches a maximum $z$-value at $(r,z)=(2.05, 0.92172)$.

\begin{figure}
  \centering
  \includegraphics[scale=.8]{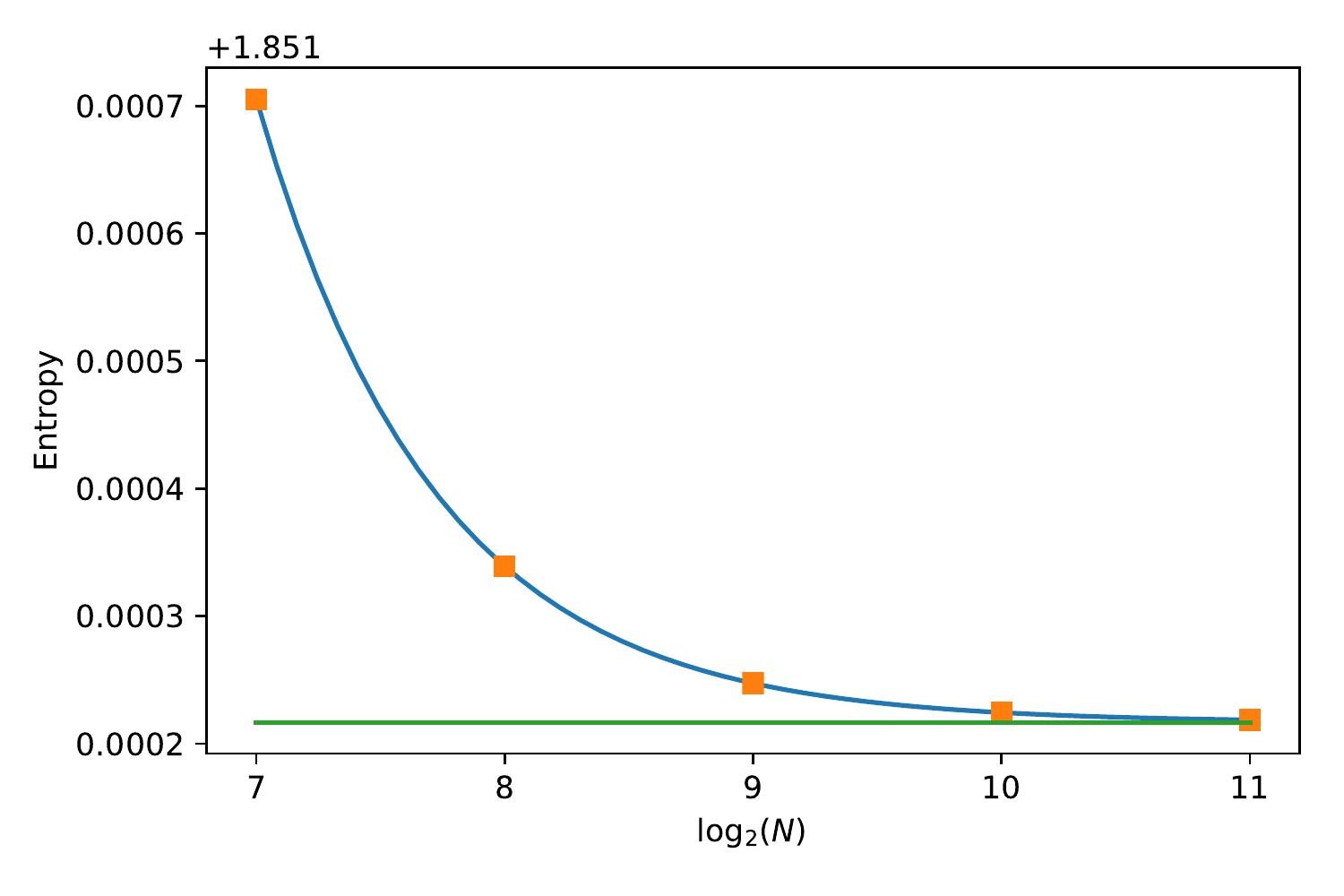}
  \caption{The entropy of the Angenent torus, computed using $128$, $256$, $512$, $1024$, and $2048$ points. The values appear to lie on an exponential curve converging to $1.8512167$.}
  \label{fig:Fvalues}
  \vspace\floatsep
  \includegraphics[scale=.8]{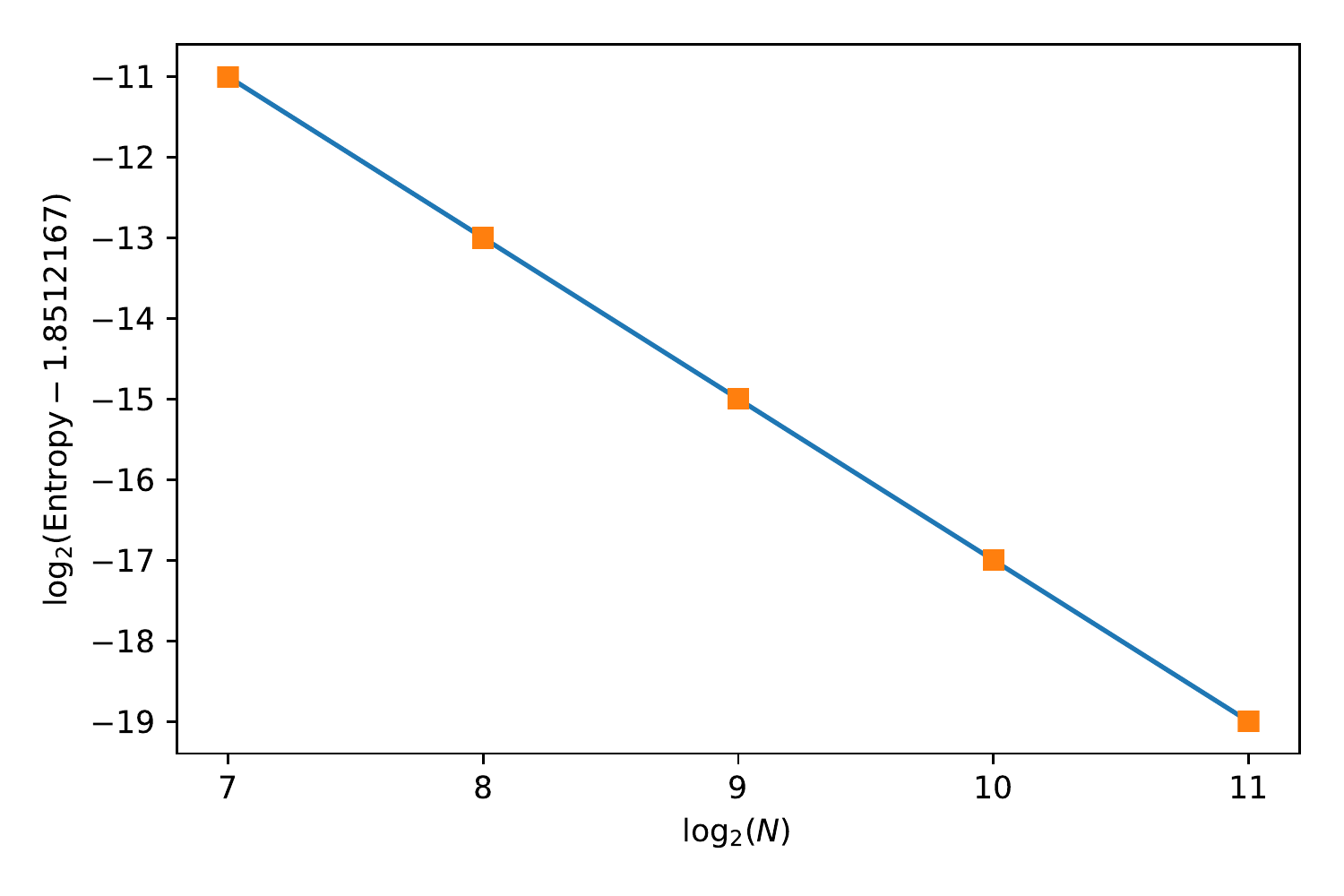}
  \caption{A log-log plot of the estimated error in the entropy value versus the number of points on the trajectory. The slope of the fit line is $-1.998$, suggesting a quadratic rate of convergence.}
  \label{fig:FvaluesLog}
\end{figure}

For each of these values of the number of points $N$, we used the closed discrete geodesic we found to estimate the entropy of the Angenent torus. The results are plotted in Figure \ref{fig:Fvalues}. All five values are just above $1.851$. The observed rate of convergence suggests that our estimates for the entropy converge quadratically. That is, as the number of points doubles, the accuracy of the estimate appears to improve by a factor of four, which is consistent with \cite[Example 2.3.2]{mw01}. See Figure \ref{fig:FvaluesLog}. Our most accurate estimate for the entropy, using $2048$ points, is $1.8512186$, and this numerical error analysis suggests that it overestimates the true value by about $0.0000019$.

We also applied our method to the sphere and the cylinder, whose entropy values were computed by Stone in \cite{s94}. Using $N=256$ points, we estimated the entropy of the sphere to be $1.471528$, whereas the true value is $\frac4e\approx 1.471518$. Using $N=1024$ points, we estimated the entropy of the cylinder to be $1.5200$, whereas the true value is $\sqrt{\frac{2\pi}e}\approx1.5203$. Neither of these cross-sections are closed curves in the $(r,z)$ half-plane, so we specified boundary conditions, solving Equations \eqref{eq:open} instead of Equations \eqref{eq:closed}. For the sphere, we specified boundary conditions on the $z$-axis. For the cylinder, we specified boundary conditions at $z=\pm\infty$. Note that the metric $g$ goes to zero at the $z$-axis and at infinity. Thus, near the $z$-axis or near infinity, changes in the discrete curve have only a very small effect on the curve's length, and hence only a very small effect on the value of the discrete action. As a consequence, we saw a numerical loss of accuracy in the computed discrete geodesics as they approached the $z$-axis or infinity, but this loss of accuracy did not affect our estimates for the geodesics' length. We do not encounter these issues with the Angenent torus, as its cross-section stays away from the $z$-axis and from infinity.

To place our numerical result of $1.85$ in context, we present it alongside the entropy values of the sphere and the cylinder computed by Stone \cite{s94}.
\begin{align*}
  \lambda(\mathbb R^2)&=1,\\
  \lambda(S^2)&\approx1.47,\\
  \lambda(S^1\times\mathbb R)&\approx1.52,\\
  \lambda(\text{Angenent torus})&\approx1.85.
\end{align*}
Thus, the entropy of the Angenent torus is larger than that of the cylinder, as we expect from the fact that we can perturb the Angenent torus to a surface that develops a cylindrical singularity under mean curvature flow. We also see that the entropy of the Angenent torus is smaller than that of two planes, as we know from \cite{dn18}.

\section{Future work}\label{sec:futurework}
There are several further directions for this work.

\begin{question}
  In this paper, we consider the two-dimensional Angenent torus, but Angenent's work applies in any dimension. By adjusting the metric $g$ appropriately, a geodesic in the $(r,z)$ half-plane will give the cross-section of a self-shrinking surface $S^1\times S^{n-1}\subset\mathbb R^{n+1}$. What are the entropies of these surfaces? What value does the entropy approach as $n$ becomes large?
\end{question}

\begin{question}
  An important question in the study of self-shrinkers is their index. Generic perturbations of a self-shrinker will have larger entropy, but if we perturb it in particular ways, the entropy may decrease. The \emph{index} of a self-shrinker is the dimension of this space of perturbations that decrease the entropy. Work of Liu \cite{l16} shows that the index of the Angenent torus is at least three, but its exact value is unknown.

  This question is amenable to numerical experimentation. The Hessian of the discrete action at the critical curve is a finite-dimensional matrix, and the eigenvectors of this matrix with negative eigenvalues represent perturbations of the cross-section that decrease the entropy, once we take dilations and translations into account. These perturbations of the cross-section yield rotationally symmetric perturbations of the Angenent torus, but one can also use perturbations of the cross-section to understand non-rotationally symmetric perturbations by making use of a result of Liu \cite{l16} that says that it suffices to consider perturbations of the torus of the form $ue^{ik\theta}$, where $u$ is rotationally symmetric. 
\end{question}

\begin{question}
  In this article, we numerically assessed the accuracy of our value for the entropy by looking at its rate of convergence as the number of points on the trajectory grows. This evidence suggests that our value is within $2\times10^{-6}$ of the true value. One would like to rigorously prove such an error bound. One approach would be to first bound the $C^1$ difference between the discrete Lagrangian \eqref{eq:discreteLagrangian} and the exact discrete Lagrangian from Definition \ref{def:exactDiscreteLagrangian}. At its core, this step amounts to understanding how well the distance between two nearby points with respect to $g$ is approximated by the Euclidean distance times the conformal factor. Then, using the Hessian of the discrete action, one could bound how far the critical curve and critical value can move when the discrete action is replaced by the exact discrete action.
\end{question}

\section{Acknowledgements}
I would like to thank Jacob Bernstein for proposing this problem at the Geometric Analysis Conference at Rutgers. I would also like to thank Bill Minicozzi and Ari Stern for their comments on this work.

\bibliographystyle{plain}
\bibliography{meancurvature}

\end{document}